\documentclass[11pt,a4paper]{amsart}
\usepackage{graphicx}
\usepackage{url}
\usepackage[latin1]{inputenc}
\usepackage[T1]{fontenc}
\usepackage[english]{babel}
\usepackage{mathscinet}
\newcommand*{\A}{\mathrm{A}'}
\newcommand*{\Ap}{\mathrm{A}^{\mathrm{p}}}
\newcommand*{\Aph}{\mathrm{A}^{\mathrm{ph}}}
\newcommand*{\Ar}{\mathcal{A}}
\newcommand*{\C}{\mathbb{C}}
\newcommand*{\Cr}{\mathcal{C}}
\def\dd{\mathop{\mathrm{d}\null}%
       \mskip-\thinmuskip\mathord{\null}}
\newcommand*{\E}{\mathrm{E}}
\newcommand*{\Ep}{\mathrm{E}^{\mathrm{p}}}
\newcommand*{\Eph}{\mathrm{E}^{\mathrm{ph}}}
\newcommand*{\fm}[2]{\mathrm{M}_{#1}[\Gamma_0(#2)]}
\DeclareMathOperator\GL{GL}
\newcommand*{\HD}{\mathcal{H}(2)}
\DeclareMathOperator{\id}{id}
\let\plainIm = \Im
\def\Im{\mathop{\plainIm\mkern-2mu \mathit{m}}\nolimits}

\DeclareMathOperator{\lcm}{lcm}
\DeclareMathOperator{\Per}{Per}
\newcommand*{\ppmod}[1]{~(\mathrm{mod}\, #1)}
\newcommand*{\pk}{\mathcal{H}}
\newcommand*{\Q}{\mathbb{Q}}
\newcommand*{\qm}[3]{\widetilde{\mathrm{M}}_{#1}%
  [\Gamma_0(#2)]^{\leq  #3}}
\newcommand*{\R}{\mathbb{R}}
\newcommand*{\rquotient}[2]{\left.\raisebox{0.2ex}{$#1$}%
   /\raisebox{-0.2ex}{$#2$}\right.}
\DeclareMathOperator\SL{SL}
\newcommand*\T{\mathbf{T}}
\newcommand*{\trsf}[3]{(#1\underset{#2}{\vert}#3)}
\newcommand*{\un}{\mathbb{1}}
\newcommand*{\Z}{\mathbb{Z}}
\newcommand*{\Zp}{\Z_{>0}}
\newcommand*{\Zpn}{\Z_{\geq 0}}
\theoremstyle{plain}
\newtheorem{prop}{Proposition}
\newtheorem{thm}[prop]{Theorem}
\newtheorem{dfn}[prop]{Definition}
\newtheorem{lem}[prop]{Lemma}
\newtheorem{cor}[prop]{Corollary}
\theoremstyle{remark}
\newtheorem*{rem}{Remark}
\title{Orbitwise countings in $\HD$ and quasimodular forms}
\author{Samuel Leli\`evre}
\address{%
Samuel Leli\`evre,
Mathematics Institute,
University of Warwick,
Coventry CV4 7AL,
United Kingdom
}
\email{samuel.lelievre@polytechnique.org}
\urladdr{http://carva.org/samuel.lelievre/}
\author{Emmanuel Royer}
\address{%
Emmanuel Royer,
I3M, UMR CNRS 5149, Universit\'e Montpellier~2, case 51,
Place Eug\`ene Bataillon, F-34095 Montpellier cedex 5, France%
}
\email{royer@math.univ-montp2.fr}
\address{%
Universit\'e Montpellier III,
MIAp,
F-34199 Montpellier cedex,
France
}
\email{emmanuel.royer@polytechnique.org}
\urladdr{http://carva.org/emmanuel.royer/}
\date{April 2006}
\begin{document}

\begin{abstract}
We prove formulae for the countings by orbit of square-tiled surfaces of genus
two with one singularity.  These formulae were conjectured by
Hubert \& Leli\`evre. We show that these countings admit quasimodular 
forms as generating functions.
\end{abstract}


\maketitle
\thispagestyle{empty}

\footnotesize
\setcounter{tocdepth}{1}
\tableofcontents
\normalsize

\section{Introduction}

The main result of this paper is the proof of a conjecture of Hubert
\& Leli\`evre.

\begin{thm}\label{thm:pal}
For odd $n$, the countings by orbit of primitive square-tiled surfaces
of the stratum $\HD$ tiled with $n$ squares are the following.
Orbit $A_n$ contains
\[
 a_n^{\mathrm{p}}= \frac{3}{16}(n-1)n^2\prod_{p\mid n}\bigl(1-\frac{1}{p^2}\bigr)
\]
primitive surfaces with $n$ squares and orbit $B_n$ contains
\[
b_n^{\mathrm{p}}=\frac{3}{16}(n-3)n^2\prod_{p\mid n}\bigl(1-\frac{1}{p^2}\bigr)
\]
primitive surfaces with $n\geq 3$ squares and $0$ for $n=1$ square.
\end{thm}
\begin{rem}
The notation $\prod_{p\mid n}$ indicates a product over prime divisors of $n$.
The superscript $\mathrm{p}$ is here to emphasize primitivity.
\end{rem}

Theorem~\ref{thm:pal} can also be expressed in terms of
quasimodularity of the generating functions of the countings. More
precisely:
\begin{cor}\label{cor:qmodat}
For any odd positive integer $n$, the number $a_n$ of $n$-square-tiled surfaces of type A in $\HD$, primitive or not, 
is the $n$th coefficient of the quasimodular form
\[
\sum_{n=0}^{+\infty}a_n\exp(2i\pi nz)
=
\frac{1}{1280}\left[E_4(z)+10\frac{\dd}{2i\pi\dd z}E_2(z)\right]
\]
of weight $4$, depth $2$ on $\SL(2,\mathbb{Z})$.
\end{cor}
\begin{rem}
Functions $E_2$ and $E_4$ are the usual Eisenstein series of weight
$2$ and $4$ respectively. They are precisely defined in equations
\eqref{eq:E2} and \eqref{eq:E4}.
\end{rem}
Since the coefficients $a_n$ have no geometric meaning for even $n$,
it makes sense to consider only the odd part of the Fourier series.
Considering the odd part is the same as considering the Fourier series
twisted by a Dirichlet character of modulus $2$ (see
section~\ref{sec:twist}).  It is then natural to expect that,
similarly to the case of modular forms (see \cite[Theorem
7.4]{Iwa97}), the odd part of the Fourier series is a quasimodular
form on the congruence subgroup $\Gamma_0(4)$.  Actually, we will
prove this is the case.  Let $\Phi_2$ and $\Phi_4$ be the two modular
forms of respective weights $2$ and $4$, defined on $\Gamma_0(4)$ as
in \eqref{eq:lesphi} and \eqref{eq:losphi}.
\begin{thm}\label{thm:labonneserie}
The Fourier series
\[
\sum_{n\in 2\Zpn+1}a_n\exp(2i\pi nz)
\] 
is the quasimodular form of weight $4$ and depth $1$ on $\Gamma_0(4)$
defined by
\[
\frac{1}{1280}
\left[E_4(z)
-
9E_4(2z)
+8E_4(4z)
-15\frac{\dd}{2i\pi\dd z}\Phi_2(z)
+15\frac{\dd}{2i\pi\dd z}\Phi_4(z)\right]
.
\]
\end{thm}
\begin{rem}
This theorem will be proved in section~\ref{sec:labonneserie}.  It is
interesting to note that forgetting the artificial terms of even order
results in a lesser depth, that is, in a simplified
modular situation. (A modular form is a quasimodular form of
depth $0$ so the depth may be seen as a measure of complexity.)
\end{rem}

\begin{table}[!h]
\centering
\begin{tabular}{|c|c|c|c|c|c|c|c|c|c|c|c|c|}
\hline
$n$ & $5$ & $7$ & $9$ & $11$ & $13$ & $15$ & $17$ & $19$
& $21$ & $23$ & $25$ & $27$\\
\hline
$a^{\mathrm{p}}_n$ & 
$18$ & $54$ & $108$ & $225$ & $378$ & $504$ & $864$ & $1215$ & 
$1440$ & $2178$ & $2700$ & $3159$\\ 
\hline
$a_n$ & 
$18$ & $54$ & $120$ & $225$ & $378$ & $594$ & $864$ & $1215$ & 
$1680$ & $2178$ & $2808$ & $3630$\\ 
\hline
\end{tabular}
\caption{Number of surfaces of type A.}
\label{tab:numb}
\end{table}


Our results may be interpreted in terms of counting genus $2$ covers of the torus $\T=\rquotient{\C}{\Z+i\Z}$ with one double ramification
point (see \S\ref{sec:geba}). The general problem of counting covers with fixed ramification type of a given Riemann surface
was posed in 1891 by Hurwitz who precisely counted the covers of the sphere. In 1995, Dijkgraaf \cite{MR1363055} computed the 
generating series of the countings of degree $n$ and genus $g$ covers of $\T$ with simple ramification over distinct points, 
weighted by the inverse of the number of automorphisms. Kaneko \& Zagier \cite{KZ} introduced the notion of quasimodular forms
and proved that the generating series computed by Dijkgraaf was quasimodular of weight $6g-6$ on $\SL(2,\Z)$. 
The case of arbitrary ramification over a single point was studied by Bloch \& Okounkov \cite{MR1742353}. 
They proved that the countings lead to linear combinations of quasimodular forms of weight less than or equal to $6g-6$.
This was used by Eskin \& Okounkov \cite{EO01} to compute volumes of the strata of moduli spaces of translation surfaces (see
also \cite{Zo}). The $\SL(2,\Z)$ orbits of square-tiled surfaces were studied by Hubert \& Leli\`evre in the case of a prime
number of squares \cite{HL1} and by McMullen \cite{Mc} in the general case.

Up to a multiplicative constant factor, our counting functions are the orbifold Euler characteristics of Teichm{\"u}ller curves.
Matt Bainbridge independently obtained results similar to ours in this setting. 

The moduli space of holomorphic 1-forms on complex curves of a
fixed genus $g$ can be considered as a family of flat structures of a special
type on a surface of genus $g$. The group $\GL(2,\R)$ acts naturally on
the moduli space; its orbits, called Teichm{\"u}ller discs, project 
to the moduli space of curves as complex geodesics for
the Teichm{\"u}ller metric. A typical flat surface has no symmetry; its
stabilizer in $\GL(2,\R)$ is trivial; the corresponding Teichm{\"u}ller disc is
dense in the moduli space. For some flat surfaces (called Veech surfaces)
the stabiliser is big (a lattice) so that the corresponding Teichm{\"u}ller
disc is closed. Projections of such Teichm{\"u}ller discs, called Teichm{\"u}ller curves, play the
role of ``closed complex geodesics''.

The main lines of the proof of theorem~\ref{thm:pal} are the
following.  In section~\ref{sec:suty}, we evaluate the number
$a_n^{\mathrm{p}}$ in terms of sums over sets defined by complicated
arithmetic conditions.  In section~\ref{sec:hlp}, we relate these
coefficients $a_n^{\mathrm{p}}$ to sums of sums of divisors of the
form
\[
\sum_{\substack{(a,b)\in\Zp^2\\ka+b=n}}
\sigma_1(a)\sigma_1(b).
\]
For the computation of these sums, we use, in section~\ref{sec:qmodu},
the notion of quasimodular forms on congruence subgroups (introduced
by Kaneko \& Zagier in \cite{KZ}) and we take advantage of the fact
that the spaces of quasimodular forms have finite dimension to
linearise the above sums.  Here, linearising means expressing them as
linear combinations of sums of powers of divisors.  Having obtained a
series whose odd coefficients are the numbers $a_n$, we introduce the
notion of twist of a quasimodular form by a Dirichlet character, to
construct a new quasimodular form generating series without artificial
Fourier coefficients.

\noindent\textbf{Thanks} A first version of this text was written as both authors were visiting the Centre 
de Recherches Math\'ematiques de Montr\'eal, and circulated under the title ``Counting integer points by Teichm{\"u}ller 
discs in $\Omega\mathcal{M}(2)$''. 
This version was written as the second author was visiting the Warwick Mathematics Institute.
We thank both institutions for the good working conditions. 
\section{Geometric background}\label{sec:geba}

\subsection{Square-tiled surfaces}
A square-tiled surface is a collection of unit squares endowed with 
identifications of opposite sides: each top side is identified to a 
bottom side and each right side is identified to a left side. In 
addition, the resulting surface is required to be connected.
A square-tiled surface tiled by $n$ squares is also a degree $n$
(connected) branched cover of the standard torus
$\rquotient{\C}{\Z+i\Z}$ with a single branch point.

Given a square-tiled surface, to each vertex can be associated an 
angle which is a multiple of $2\pi$ (four or a multiple of four 
squares can abutt at each vertex). If $(k_i+1)2\pi$ is the angle at vertex $i$, the Gauss--Bonnet 
formula implies that
\[
\sum_{i=1}^sk_i=2g-2
\]
where $g$ is the genus of the surface and $s$ the total number of vertices.

\begin{figure}[ht!]
\centering
\includegraphics[width=0.5\linewidth]{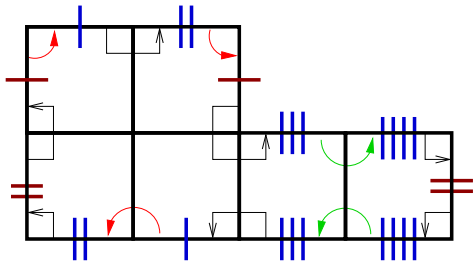}
\caption{Surface with one angle of $6\pi$}
\label{fig:Angles}
\end{figure}

Surfaces can be sorted according to strata
$\mathcal{H}(k_1,\dotsc,k_s)$. Square-tiled surfaces are integer points of these strata. In this paper we are concerned with
surfaces in $\HD$, that is, with a single ramification point of angle
$6\pi$. A surface tiled by $n$ squares in $\HD$ is a degree $n$ branched cover
of the torus $\rquotient{\C}{\Z+i\Z}$ with one double ramification
point.

\subsection{Cylinder decompositions}
Given any square-tiled surface, each horizontal line on the surface through the interior of a square is closed, and neighbouring 
horizontal lines are also closed.  Thus closed horizontal lines come in families forming
cylinders and the surface decomposes into such cylinders bounded by
horizontal saddle connections (segments joining conical
singularities).

\begin{figure}[ht!]
\centering
\includegraphics[width=0.8\linewidth]{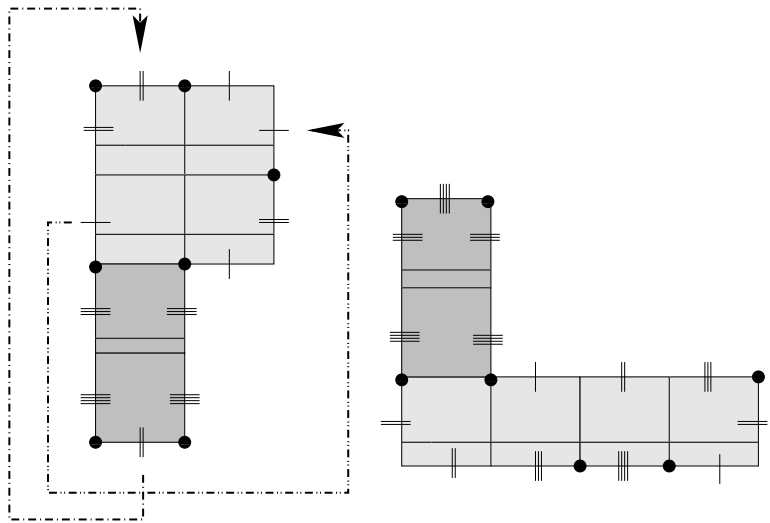}
\caption{Decomposition of a surface into two cylinders}
\label{fig:deccyl}
\end{figure}

Here we explain how to enumerate square-tiled surfaces in $\HD$ with a
given number of squares, by giving a system of coordinates for them.
We include this discussion for the sake of completeness, although 
these coordinates have already been used in \cite{Zo}, \cite{EMS}, \cite{HL1}.

We represent surfaces according to their cylinder decompositions.
Cylinders of a square-tiled surface are naturally represented as rectangles. One can cut a triangle
from one side of such a rectangle and glue it back on the other side according to the identifications to produce
a parallelogram with a pair of horizontal
sides (each made of one or several saddle connections), and a pair of
identified nonhorizontal parallel sides. A square-tiled surface in $\HD$ has one or two cylinders \cite{Zo} and can always
be represented as in figure \ref{fig:twun} or \ref{fig:twdeux}.
Each cylinder has a height and a width and in
addition a twist parameter corresponding to the possibility of rotating the saddle connections of the top or
bottom of the cylinders before performing the identifications. 

\subsection{One-cylinder surfaces}

\begin{figure}[ht!]
\centering
\includegraphics[width=\linewidth]{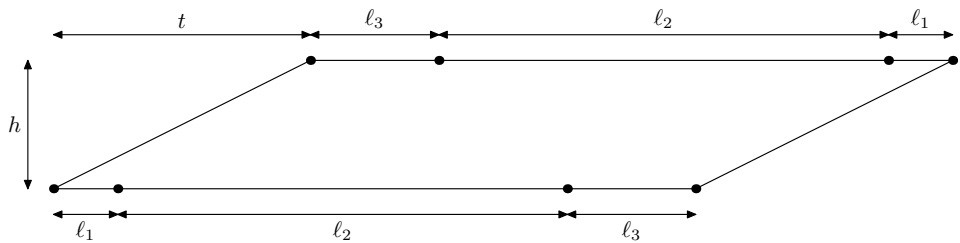}
\caption{One-cylinder surface}
\label{fig:twun}
\end{figure}

For one-cylinder surfaces in $\HD$, we have on the bottom of the
cylinder three horizontal saddle connections, and the same saddle
connections appear on the top of the cylinder in reverse order; we
denote by $\ell$ the width of the cylinder and $\ell_1$, $\ell_2$,
$\ell_3$ the lengths of the saddle connections, numbered so that they
appear in that order on the bottom side and in reverse order on the top 
side. See figure~\ref{fig:twun}.

For each choice of $(\ell_1,\ell_2,\ell_3)$ with $\ell_1+\ell_2+\ell_3=n$, if $\ell_1$, $\ell_2$, $\ell_3$ are
not all equal, there are $\ell$ possible values of the twist $t$ giving different surfaces. But, the three possible 
cyclic permutations of $(\ell_1,\ell_2,\ell_3)$ yields the same set of surfaces. So, to make coordinates uniquely
defined, we require that $(\ell_1,\ell_2,\ell_3)$ has least lexicographic order among its cyclic permutations.
For countings, it is simpler to ignore this point and to divide by $3$ at the end. 

If $\ell_1$, $\ell_2$, $\ell_3$ are all equal (and thus worth 
$\ell/3$), there is only one cyclic permutation of
$(\ell_1,\ell_2,\ell_3)$ but only $\ell/3$ values of the twist $t$ give 
different surfaces.

The parameters we have used satisfy:
\begin{align*}
\ell &\mid n\\
\ell_1+\ell_2+\ell_3&=\ell\\
0\leq t < \ell\text{ or }\ell/3.
\end{align*}

\begin{rem}
From this description of coordinates, we conclude that the number of
one-cylinder surfaces in $\HD$ tiled with $n$ squares is (see \cite{EMS})
\[
\frac{1}{3}\sum_{\ell\mid
  n}\sum_{\substack{(\ell_1,\ell_2,\ell_3)\in\Zp^{3}\\
  \ell_1+\ell_2+\ell_3=\ell}}\ell.
\]
\end{rem}

\subsection{Two-cylinder surfaces}

\begin{figure}[ht!]
\centering
\includegraphics[width=\linewidth]{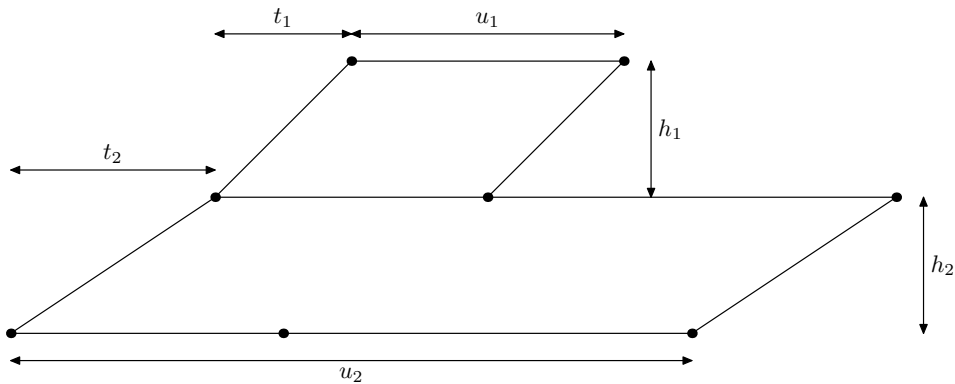}
\caption{Two-cylinder surface}
\label{fig:twdeux}
\end{figure}

Given a two-cylinder surface in $\HD$, one of its cylinders (call it cylinder
$1$) has one saddle connection on the top and one saddle connection
(of same length) on the bottom, while the other one (call it cylinder
$2$) is bounded by two saddle connections on the top and two saddle
connections on the bottom. See figure~\ref{fig:twdeux}.

For each of cylinders $1$ and $2$, there are three parameters: the 
height $h_i$, the width $u_i$, the twist $t_i$.

Given two heights $h_1$ and $h_2$, two widths $u_1<u_2$, and two
twists $t_1$, $t_2$ with $0 \leq t_i < u_i$, there exists a unique 
surface in $\HD$ with two cylinders having 
$(h_1,h_2,u_1,u_2,t_1,t_2)$ as parameters. The number of squares is 
then $h_1u_1 + h_2u_2$.

\begin{rem}
From this system of coordinates one deduces (see \cite{EMS}) that the 
number of two-cylinder surfaces in $\HD$ tiled by $n$ squares is
\[
\sum_{\substack{(h_1,h_2,u_1,u_2)\in\Zp^{4}\\
u_1<u_2\\ h_1u_1+h_2u_2=n}}u_1u_2.
\]
\end{rem}

\subsection{Lattice of periods}

The lattice of periods of a square-tiled surface is the rank two 
sublattice of $\Z^2$ generated by its saddle connections.

\begin{lem}
A square-tiled surface is translation-tiled by a parallelogram 
if and only if this parallelogram is a fundamental domain for
a lattice containing the surface's lattice of periods.
\end{lem}

\begin{proof}
Decompose the surface into polygons with vertices at the conical
singularities. The sides of these polygons are saddle connections
and together generate the lattice of periods. The tiling of the
plane by parallelograms which are a fundamental domain for this
lattice (or any rank two lattice of the plane containing it) 
yields a tiling of the translation surface by such parallelograms.
\end{proof}

\begin{rem}
The previous lemma implies that the area of the lattice of periods
divides the area of the surface it comes from.
\end{rem}

We will describe the basis of the lattice of periods by use of the following lemma \cite[Chapitre 7]{Ser97}.
\begin{lem}\label{lem:base}
Let $\Lambda$ be a sublattice of $\Z+i\Z$ of index $d$. Then there exists a unique triple of integers $(a,t,h)$ with
$a\geq 1$, $ah=d$ and $0\leq t\leq a-1$ such that
\[
\Lambda=(a,0)\Z\oplus (t,h)\Z.
\]
\end{lem}

\begin{rem}
Let $S$ be a square-tiled surface in $\HD$, and let
$\bigr(\begin{smallmatrix}a&t\\0&h\\\end{smallmatrix}\bigl)$
be the matrix corresponding to its lattice of periods.
If $S$ is one cylinder then $h$ is the height of its unique cylinder ; if $S$ is two-cylinder, with cylinders of height $h_1$ 
and $h_2$, then $h = (h_1, h_2)$.
\end{rem}

\begin{dfn}
A square-tiled surface is called \emph{primitive} if its lattice of
periods is $\Z^2$, in other words if
$\bigr(\begin{smallmatrix}a&t\\0&h\\\end{smallmatrix}\bigl)
= \bigr(\begin{smallmatrix}1&0\\0&1\\\end{smallmatrix}\bigl)$.
\end{dfn}

\begin{dfn}
A square-tiled surface is called \emph{primitive in height} if $h=1$.
\end{dfn}

The linear action of $\GL(2,\Q)^{+}$ on $\R^2$ induces an action of 
$\GL(2,\Q)^{+}$ on square-tiled surfaces. This action preserves 
orientation. The action of $\SL(2,\Z)$ preserves the number of square 
tiles, and preserves primitivity.

Hubert \& Leli\`evre have shown that if $n\geq 5$ is prime, the
square-tiled surfaces in $\HD$ tiled with $n$ squares (necessarily
primitive since $n$ is prime) form into two orbits under $\SL(2,\Z)$,
denoted by $A_n$ and $B_n$.

If $n$ is not prime and $n\geq6$, not all surfaces tiled by $n$
squares are primitive, and if $n$ has many divisors these surfaces
split into many orbits under $\SL(2,\Z)$, most of them lying in orbits
under $GL(2,\Q)^{+}$ of primitive square-tiled surfaces with fewer
squares.  There can be an arbitrary number of such ``artificial''
$\SL(2,Z)$-orbits.  Artificial orbits consist only of nonprimitive square-tiled
surfaces, since the action of $\SL(2,\Z)$ preserves primitivity.

Let $n$ be an odd integer. We can distinguish two types of surfaces 
among surfaces tiled by $n$ squares in $\HD$. These two types are 
distinguished by Weierstrass points, as follows (see \cite{HL1}).

On a surface in $\HD$, the matrix 
$\bigl(\begin{smallmatrix}-1&0\\0&-1\\\end{smallmatrix}\bigr)$
induces an involution which can be shown to have six fix points, called the Weierstrass points of the surface.
It is easy to show that for a square-tiled surface these points have coordinates in 
$\frac{1}{2}\Z$.
The type invariant is determined by the number of Weierstrass points 
which have both their coordinates integer:
\begin{itemize}
\item
a surface is of type A if it has one integer Weierstrass point;
\item
a surface is of type B if it has three integer Weierstrass points.
\end{itemize}

\begin{rem}
We give an interpretation in terms of orbits.
Consider the orbit under $\GL(2,\Q)^{+}$ of a surface $S$ tiled by 
$n$ squares. Then
\begin{itemize}
\item
the primitive square-tiled surfaces in this orbit all have the same 
number of squares, say $d$,
\item
the action of $\GL(2,\Q)^{+}$ restricts to an action of $\SL(2,\Z)$ 
on these primitive square-tiled surfaces;
\item
these primitive surfaces form an orbit under $\SL(2,\Z)$;
\item
McMullen extended the result of Hubert \& Leli\`evre by showing that
if $d\geq5$ is odd, the set of primitive square-tiled surfaces in
$\HD$ tiled with $d$ squares is partitioned in two orbits under 
$\SL(2,\Z)$ denoted $A_d$ and $B_d$;
\item
the type can be read on these primitive square-tiled surfaces.
\end{itemize}

\end{rem}
 
\section{Sum-type formulae for the orbitwise countings}\label{sec:suty}

\subsection{From primitive to non primitive countings}

In this section, we establish relations between countings of primitive surfaces, countings of height primitive
surfaces and countings of (non necessarily primitive) surfaces. For any integer $\ell$, the function $\sigma_\ell$ is defined by
\begin{equation}\label{eq:defsig}
\sigma_\ell(n)
=
\begin{cases}
\sum_{d\mid n}d^\ell & \text{if $n\in\Zp$}\\
0 & \text{otherwise.}
\end{cases}
\end{equation}

For $n\in\Zp$, we define $\E_n$ as the set of $n$ squares surfaces in $\HD$, $\Ep_n$ as its subset of primitive surfaces
and $\Eph_n$ as its subset of primitive in height surfaces.
For $d\in\Zp$, we note $\Lambda_d$ the set of sublattices of $\Z+i\Z$ of index $d$.
The description of surfaces by primitive surfaces is given by the following lemma.
\begin{lem}\label{lem:EpE}
For $n\in\Zp$, we have the following bijection
\[
\E_n\simeq\bigcup_{d\mid n}\Ep_{n/d}\times\Lambda_d.
\]
\end{lem}
\begin{proof}
Let $S\in\E_n$ and $d$ be the index in $\Z+i\Z$ of its lattice of periods $\Per(S)$. Then $d\mid n$ and $\Per(S)\in\Lambda_d$.
With the notations of lemma \ref{lem:base}, we write $\Per(S)=(a,0)\Z\oplus (t,h)\Z$. To $S$ we associate a surface tiled by
$n/d$ squares:
\[
S'=\begin{pmatrix}a & t\\ 0 & h\end{pmatrix}^{-1}S.
\]
The lattice of periods of $S'$ is $\Z+i\Z$ so that it is primitive. Conversely, let $S'\in\Ep_{n/d}$ and $\Lambda\in\Lambda_d$.
With the notations of lemma \ref{lem:base}, we write $\Lambda=(a,0)\Z\oplus (t,h)\Z$. Then
\[
S=\begin{pmatrix}a & t\\ 0 & h\end{pmatrix}S'
\]
has $n=ah$ squares.
\end{proof}
\begin{cor}
For $n\in\Zp$, we have
\[
\#\E_n=\sum_{d\mid n}\sigma_1(d)\#\Ep_{n/d}.
\]
\end{cor}
\begin{proof}
By lemma \ref{lem:base}, we have
\[
\#\Lambda_d=\sum_{\substack{(a,t,h)\in\Zpn^3\\ ah=d\\ 0\leq t<a}}1=\sigma_1(d).
\]
\end{proof}
We recall that a surface $S$ is primitive in height if $h=1$ with the notations of lemma \ref{lem:base}. That is, its
lattice of periods is $\Per(S)=(a,0)\Z+(t,1)\Z$ with $a\geq 1$ and $0\leq t\leq a-1$. We write $\Lambda'_d$ for the set
of these lattices having index $d$ (implying $d=a$). We have $\#\Lambda'_d=d$. Similary to lemma \ref{lem:EpE} we have
\begin{lem}\label{lem:EphE}
For $n\in\Zp$, we have the following bijection
\[
\E_n\simeq\bigcup_{d\mid n}\Eph_{n/d}\times\Lambda'_d.
\]
\end{lem}
\begin{cor}
For $n\in\Zp$, we have
\[
\#\E_n=\sum_{d\mid n}d\cdot\#\Eph_n.
\]
\end{cor}

We deduce the same result for surfaces of type A. For odd $n$, define $\A_n$ as the set of $n$-square-tiled surfaces of type A in $\HD$,
$\Ap_n$ as its subset of primitive surfaces (which coincides with the $\SL(2,\Z)$-orbit $A_n$) and $\Aph_n$ as its subset of 
height-primitive surfaces.
\begin{lem}\label{lem:ApA}
For $n\in\Z$ odd, we have the following bijection
\[
\A_n\simeq\bigcup_{d\mid n}\Ap_{n/d}\times\Lambda_d.
\]
\end{lem}
\begin{proof}
We recall that the type of a surface is characterized by the number of its Weierstrass points with integer coordinates. To deduce 
lemma \ref{lem:ApA} from lemma \ref{lem:EpE} it then suffices to prove that a Weierstrass point $P$ has half-integer
coordinates\footnote{Meaning in $\frac{1}{2}\Z^2$ but not in $\Z^2$} in a basis determined by $\Per(S)$ if and only if its image by 
the bijection of lemma \ref{lem:EpE} has half-integer coordinates
in the canonical basis of $\Z+i\Z$. Let $S\in\E_n$, $\Per(S)=(a,0)\Z\oplus (t,h)\Z$ its lattice of periods with the notations of lemma
\ref{lem:base}. We set
\[
M=\begin{pmatrix}a & t\\ 0 & h\end{pmatrix}.
\]
Let $P$ a Weierstrass point in $S$, we assume that its coordinates in the basis of $\Per(S)$ are $(\ell/2,m/2)$ with $m$ and $n$
not simultaneously even. The coordinates of $P$ in $\Z+i\Z$ are therefore $(a\ell+mt,mh)/2$, hence, those of $M^{-1}P$ in
$M^{-1}S$ are $(\ell/2,m/2)$ in the standard basis of $\Z+i\Z$. 
\end{proof}
\begin{cor}
For $n\in\Zp$, we have
\[
a_n=\sum_{d\mid n}\sigma_1(d)a_{n/d}^{\mathrm{p}}.
\]
\end{cor}
\begin{lem}\label{lem:AphA}
For $n\in\Z$ odd, we have the following bijection
\[
\A_n\simeq\bigcup_{d\mid n}\Aph_{n/d}\times\Lambda'_d.
\]
\end{lem}
\begin{cor}
For $n\in\Zp$, we have
\[
a_n=\sum_{d\mid n}da_{n/d}^{\mathrm{ph}}.
\]
\end{cor}

To express the number of primitive surfaces in terms of the numbers of primitive in height ones, we recall some basic facts 
on $L$-functions. For an arithmetic function $f$, we define
\[
L(f,s)=\sum_{n=1}^{+\infty}f(n)n^{-s}.
\]
If $\id$ denotes the identity function, we have
\[
L(\id^\ell f,s)=L(f,s-\ell).
\]
For $f$ and $g$ two arithmetic functions with convolution product $f\ast g$, we have
\[
L(f\ast g,s)=L(f,s)L(g,s).
\] 
The constant equal to $1$ function is denoted by $\un$ and, we have
\[
L(\un,s)=\zeta(s)
\]
the Riemann $\zeta$ function. Moreover
\[
L(\mu,s)=\frac{1}{\zeta(s)} \text{ and }
L(\sigma_k,s)=\zeta(s)\zeta(s-k).
\]

\begin{lem}\label{lem:prph}
Let $n\in\Zp$. Then
\[
a_n^{\mathrm{p}}=\sum_{d\mid n}\mu(d)a_{n/d}^{\mathrm{ph}}.
\]
\end{lem}
\begin{proof}
Lemma \ref{lem:AphA} is then
\[
L(a_n,s)=\zeta(s-1)L(a^{\mathrm{ph}},s)
\]
and lemma \ref{lem:ApA} is
\[
L(a_n,s)=\zeta(s)\zeta(s-1)L(a^{\mathrm{p}},s).
\]
We deduce 
\[
L(a^{\mathrm{p}},s)=\frac{1}{\zeta(s)}L(a^{\mathrm{ph}},s)
\]
hence the result.
\end{proof}

Next, we give sum-type formulae for the number of surfaces in $\Aph_n$.
\begin{prop}\label{prop:countUn}
Let $n\in\Zp$, the number of height-primitive one-cylinder surfaces with $n$ squares in $\HD$ of type A is
\[
\frac{1}{3}\sum_{\substack{\ell_1,\ell_2,\ell_3 \text{ odd}\\\ell_1+\ell_2+\ell_3=n}}n.
\]
\end{prop}
\begin{proof}
See figure~\ref{fig:WeUnCyl}.
Since the cylinder is primitive in height, it has height $1$.
As proved in \cite[\S 5.1.1]{HL1}, the Weierstrass points are
\begin{itemize}
\item the saddle point, which has integer coordinates
\item two points lying on the core of the cylinder, which do not have integer coordinates
\item the midpoints of the three saddle connections, each of these points having integer coordinates if and only if
the corresponding saddle connection has even length.
\end{itemize}
\begin{figure}[ht!]
\centering
\includegraphics[width=0.7\linewidth]{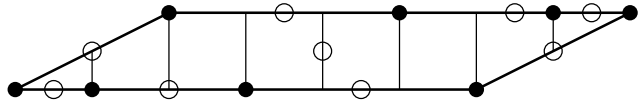}
\caption{Weierstrass points of a one-cylinder surface}
\label{fig:WeUnCyl}
\end{figure}
\end{proof}
\begin{prop}\label{prop:countDeux}
Let $n\in\Zp$, the number of height-primitive two-cylinder surfaces with $n$ squares in $\HD$ of type A is
\[\sum_{\substack{h_1,h_2,u_1,u_2\in\Zpn \\ h_1u_1 + h_2u_2 = n \\
(h_1,h_2)=1\\ h_1,h_2\text{ odd}\\ u_1<u_2}}u_1u_2%
+ 
\frac{1}{2}
\sum_{\substack{h_1,h_2,u_1,u_2\in\Zpn \\ h_1u_1 + h_2u_2 = n \\
(h_1,h_2)=1 \\ h_1\not\equiv h_2 \ppmod{2} \\ u_1<u_2 \\ \text{$u_1u_2$ even}}}u_1u_2.
\]
\end{prop}
\begin{proof} Among height-primitive two-cylinder
surfaces with parameters $h_1$, $h_2$, $u_1$, $u_2$, $t_1$, $t_2$,
such that $h_1u_1 + h_2u_2 = n$ (odd):
\begin{itemize}
\item all surfaces with $h_1$ and $h_2$ odd are of type A; 
\item all surfaces with $u_1$ and $u_2$ odd are of type B; 
\item exactly half of the remaining surfaces (with different parity for
$u_1$ and $u_2$ and for $h_1$ and $h_2$) are of type A, and half are
of type B;
\end{itemize}
for each $(h_1,h_2,u_1,u_2)$, there are $u_1 u_2$ possible twists 
(if $n$ is not prime some values of the twists may yield non primitive surfaces which is why we only require
height-primitivity).
In the case of different parities for
$h_1$ and $h_2$ and for $u_1$ and $u_2$, the product $u_1u_2$ is even 
and exactly half of the possible twists corresponds to each type.
The height-primitivity condition is just that the heights of the
cylinders have greatest common divisor equal to one.
\begin{figure}[ht!]
\centering
\includegraphics[width=0.5\linewidth]{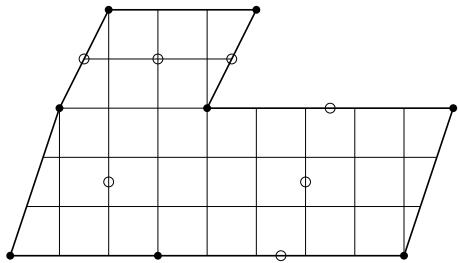}
\caption{Weierstrass points of a two-cylinder surface}
\label{fig:WeDeCyl}
\end{figure}
\end{proof}

\section{Quasimodular forms}\label{sec:qmodu}
\subsection{Motivation}

The aim of this part is the computation of sums of type
\begin{equation}\label{eq:defsk}
S_k(n)
=
\sum_{\substack{(a,b)\in\Zp^2\\ka+b=n}}
\sigma_1(a)\sigma_1(b)
\end{equation}
with $k\in\Zp$.  Here we study only the cases $k\in\{1,2,4\}$ but
the method in fact applies to every $k$ \cite{RoWil}.

Useful to the study of these sums is the weight $2$ Eisenstein series
\begin{equation}\label{eq:E2}
E_2(z)
=
1-24\sum_{n=1}^{+\infty}\sigma_1(n)e(nz)
\end{equation}
where \[e(\tau)=\exp(2i\pi\tau) \qquad (\Im\tau>0).\]
Defining
\[
H_k(z)
=
E_2(z)E_2(kz),
\]
one gets
\begin{equation}\label{eq:defhk}
H_k(z)
=
1-24\sum_{n=1}^{+\infty}\left[\sigma_1(n)+\sigma_1\left(\frac{n}{k}\right)\right]e(nz)
+576\sum_{n=1}^{+\infty}S_k(n)e(nz).
\end{equation}
We shall achieve the linearisation of $H_k$ using the theory of
quasimodular forms, developed by Kaneko \& Zagier.  The computation of $S_k(n)$
will be deduced for each $n$.

\subsection{Definition}

Let us therefore begin by surveying our prerequisites on quasimodular forms,
referring to \cite[\S 17]{MR05} for the details. 
Define
\[
\Gamma_0(N)=\left\{
\begin{pmatrix} a & b\\ c & d\end{pmatrix} \colon (a,b,c,d)\in\Z^4,\, ad-bc=1,\, N\mid c
\right\}
\]
for all integers $N\geq 1$.  In particular, $\Gamma_0(1)$ is
$\SL(2,\Z)$.  Denote by $\pk$ the Poincar\'e upper half
plane:
\[
\pk=\{ z\in\C \colon \Im z>0\}.
\]
\begin{dfn}
Let $N\in\Zp$, $k\in\Zpn$ and $s\in\Zpn$. A holomorphic function
\[
f \colon \pk \to \C
\]
is a quasimodular form of weight $k$, depth $s$ on $\Gamma_0(N)$ if
there exist holomorphic functions $f_0$, $f_1$, $\dotsc$, $f_s$ on
$\pk$ such that
\begin{equation}\label{eq:cqm}
(cz+d)^{-k}
f\left(\frac{az+b}{cz+d}\right)
=
\sum_{i=0}^sf_i(z)\left(\frac{c}{cz+d}\right)^i
\end{equation}
for all $\bigl(\begin{smallmatrix} a & b\\ c &
d\end{smallmatrix}\bigr)\in\Gamma_0(N)$ and such that $f_s$ is
holomorphic at the cusps and not identically vanishing.  By
convention, the $0$ function is a quasimodular form of depth $0$ for
each weight.
\end{dfn}
Here is what is meant by the requirement for $f_s$ to be holomorphic at
the cusps.  One can show \cite[Lemme 119]{MR05} that if $f$ satisfies
the quasimodularity condition \eqref{eq:cqm}, then $f_s$ satisfies the
modularity condition
\[
(cz+d)^{-(k-2s)}
f_s\left(\frac{az+b}{cz+d}\right)
=
f_s(z)
\]
for all $\bigl(\begin{smallmatrix} a & b\\ c &
d\end{smallmatrix}\bigr)\in\Gamma_0(N)$.  Asking $f_s$ to be
holomorphic at the cusps is asking that, for all
$M=\bigl(\begin{smallmatrix} \alpha & \beta\\ \gamma &
\delta\end{smallmatrix}\bigr)\in\Gamma_0(1)$, the function
\[
z\mapsto(\gamma z+\delta)^{-(k-2s)}f_s\left(\frac{\alpha z+\beta}{\gamma z+\delta}\right)
\]
has a Fourier expansion of the form
\[
\sum_{n=0}^{+\infty}\widehat{f}_{s,M}(n)e\left(\frac{nz}{u_M}\right)
\]
where
\[
u_M=\inf\{
u\in\Zp \colon T^u\in M^{-1}\Gamma_0(N)M
\}.
\]
In other words, $f_s$ is automatically a modular function and is
required to be more than that, a modular form of weight $k-2s$
on $\Gamma_0(N)$.  It follows that if $f$ is a quasimodular form of
weight $k$ and depth $s$, non identically vanishing, then $k$ is even
and $s\leq k/2$.

\begin{rem}
Let $\chi$ be a Dirichlet character (see \S\,\ref{sec:twist}).  If $f$
satisfies all of what is needed to be a quasimodular form except
\eqref{eq:cqm} being replaced by
\[
(cz+d)^{-k}
f\left(\frac{az+b}{cz+d}\right)
=
\chi(d)\sum_{i=0}^nf_i(z)\left(\frac{c}{cz+d}\right)^i,
\]
then, one says that $f$ is a quasimodular form of weight $k$, depth
$s$ and character $\chi$ on $\Gamma_0(N)$.
\end{rem}

The Eisenstein series $E_2$ transforms as
\[
(cz+d)^{-2}E_2\left(\frac{az+b}{cz+d}\right)
=
E_2(z)
+
\frac{6}{i\pi}\frac{c}{cz+d}
\]
under the action of any $\bigl(\begin{smallmatrix} a & b\\ c & d\end{smallmatrix}\bigr)\in\Gamma_0(1)$. Hence, $E_2$
is a quasimodular form of weight $2$ and depth $1$ on $\Gamma_0(1)$. Defining
\[
E_{N,2}(z)=E_2(Nz),
\]
one has
\[
(cz+d)^{-2}E_{N,2}\left(\frac{az+b}{cz+d}\right)
=
E_{N,2}(z)
+
\frac{6}{i\pi N}\frac{c}{cz+d}
\]
for all $\bigl(\begin{smallmatrix} a & b\\ c & d\end{smallmatrix}\bigr)\in\Gamma_0(N)$. Hence, $E_{N,2}$ is a quasimodular form
of weight $2$ and depth $1$ on $\Gamma_0(N)$. One denotes by $\qm{k}{N}{s}$ the space of quasimodular forms of weight $k$ and depth 
less than or equal to $s$ on $\Gamma_0(N)$. The space $\qm{k}{N}{0}$ is the space $\fm{k}{N}$ of modular forms of weight $k$ 
on $\Gamma_0(N)$. A recurrence on the depth implies \cite[Corollaire 121]{MR05} the equality
\begin{equation}\label{eq:baun}
\qm{k}{N}{s}
=
\bigoplus_{i=0}^s\fm{k-2i}{N}E_2^i.
\end{equation}
It is known that $\fm{2}{1}=\{0\}$. However, if $N>1$, one deduces from
\[
\C E_2\oplus\C E_{N,2}
\subset
\qm{2}{N}{1}
=
\fm{2}{N}\oplus\C E_2
\]
that $\dim\fm{2}{N}\geq 1$. By the way, for every family $(c_d)_{d\mid N}$ such that
\[
\sum_{d\mid N}\frac{c_d}{d}=0
\]
one has
\[
\left[
z\mapsto \sum_{d\mid N}c_dE_2(dz)
\right]\in\fm{2}{N}.
\]
Denote by $D$ the differential operator
\[
D=\frac{1}{2i\pi}\frac{\dd}{\dd z}.
\]
It defines a linear application from $\qm{k}{N}{s}$ to
$\qm{k+2}{N}{s+1}$.  This application is injective and strictly
increases the depth if $k>0$.  This property allows to linearise the
basis given in \eqref{eq:baun}.
\begin{lem}\label{lem:0040}
Let $k\geq 2$ even. Then,
\[
\qm{k}{N}{k/2}
=
\bigoplus_{i=0}^{k/2-1}D^i\fm{k-2i}{N}\oplus \C D^{k/2-1}E_2.
\]
\end{lem} 

\subsection{Sums of sums of divisors}

Lemma~\ref{lem:0040} allows to reach our goal by expressing the sums
$S_1$, $S_2$ and $S_4$ introduced in \eqref{eq:defsk} as follows.

\begin{prop}\label{prop:essk}
Let $n\geq 1$. Then,
\[
S_1(n)
=
\frac{5}{12}\sigma_3(n)
-\frac{n}{2}\sigma_1(n)
+\frac{1}{12}\sigma_1(n),
\]
\[
S_2(n)
=
\frac{1}{12}\sigma_3(n)
+\frac{1}{3}\sigma_3\left(\frac{n}{2}\right)
-\frac{1}{8}n\sigma_1(n)
-\frac{1}{4}n\sigma_1\left(\frac{n}{2}\right)
+\frac{1}{24}\sigma_1(n)
+\frac{1}{24}\sigma_1\left(\frac{n}{2}\right)
\]
and
\begin{multline*}
S_4(n)
=
\frac{1}{48}\sigma_3(n)
+\frac{1}{16}\sigma_3\left(\frac{n}{2}\right)
+\frac{1}{3}\sigma_3\left(\frac{n}{4}\right)
-\frac{1}{16}n\sigma_1(n)
-\frac{1}{4}n\sigma_1\left(\frac{n}{4}\right)\\
+\frac{1}{24}\sigma_1(n)
+\frac{1}{24}\sigma_1\left(\frac{n}{4}\right).
\end{multline*}
\end{prop}
\begin{proof}
We detail the proof for the expression of $S_4$.  The function $H_4$,
introduced in \eqref{eq:defhk}, is a quasimodular form of weight $4$
and depth $2$ on $\Gamma_0(4)$.  Lemma~\ref{lem:0040} gives
\[
\qm{4}{4}{2}
=
\fm{4}{4}\oplus D\fm{2}{4}
\oplus
\C DE_2.
\]
The space $\fm{4}{4}$ has dimension $3$ and contains the linearly
independent functions
\begin{align}
\label{eq:E4}
E_4(z)
&=
1+240\sum_{n=1}^{+\infty}\sigma_3(n)e(nz)\\
\notag
E_{2,4}(z)&=E_4(2z)\\
\notag
E_{4,4}(z)&=E_4(4z).
\end{align}
The space $\fm{2}{4}$ has dimension $2$ and is generated by
\begin{align}\label{eq:lesphi}
\Phi_2(z)
&=
2E_2(2z)-E_2(z)\\
\label{eq:losphi}
\Phi_4(z)
&=
\frac{4}{3}E_2(4z)-\frac{1}{3}E_2(z).
\end{align}
Hence, by the computations of the first seven Fourier coefficients, one gets
\[
H_4
=
\frac{1}{20}E_4
+
\frac{3}{20}E_{2,4}
+\frac{4}{5}E_{4,4}
+\frac{9}{2}D\Phi_4
+3DE_2.
\]
The computation of $S_4$ is then obtained by comparison of the Fourier
coefficients of this equality.  The computation of $S_2$ is obtained
along the same lines 
via the equality
\[
H_2=\frac{1}{5}E_4+\frac{4}{5}E_{2,4}+3D\Phi_2+6DE_2
\]
between forms of $\qm{4}{2}{2}$.  At last, the expression of $S_1$
is deduced from the equality
\[
E_2^2=E_4+12DE_2
\]
between forms of $\qm{4}{1}{2}$.
\end{proof}

\begin{rem}

The computation of $H_4$, which lies in the dimension $6$ vector space
with basis $\{E_4,E_{2,4},E_{4,4},D\Phi2, D\Phi_4,DE_2\}$, required
working on \emph{seven} consecutive Fourier coefficients.  We
briefly explain why, mentioning that any sequence of $6$
consecutive coefficients is not sufficient.  For any function
\[
f(z)=\sum_{n=0}^{+\infty}\widehat{f}(n)e(nz)
\] 
and any integer $i\geq 0$, define
\[
c(f,i)
=
\left(
\widehat{f}(i),
\widehat{f}(i+1),
\widehat{f}(i+2),
\widehat{f}(i+3),
\widehat{f}(i+4),
\widehat{f}(i+5)
\right),
\]
and let
\begin{align*}
v_1(i) &= c(E_4,i) &
v_2(i) &= c(E_{2,4},i) &
v_3(i) &= c(E_{4,4},i)\\
v_4(i) &= c(D\Phi_2,i)&
v_5(i) &= c(D\Phi_4,i) &
v_6(i) &= c(DE_2,i).
\end{align*}
Then, for each $i$, there exists an explicitly computable linear
relation between $v_2(i)$, $v_3(i)$, $v_4(i)$, $v_5(i)$ and $v_6(i)$.
One could think of using a basis of $\qm{4}{1}{2}$ echelonized by 
increasing powers of $e(z)$.
The same phenomenon would however appear when changing to such a basis and 
expressing the new basis elements in terms of the original basis.

\end{rem}

\subsection{Twist by a Dirichlet character}\label{sec:twist}

Recall that a Dirichlet character $\chi$ is a character of a
multiplicative group $\left(\rquotient{\Z}{q\Z}\right)^\times$
extended to a function on $\Z$ by defining
\[
\chi(n)
=
\begin{cases}
\chi\left(n\ppmod{q}\right) & \text{if $(n,q)=1$}\\
0 & \text{otherwise}
\end{cases}
\]
(see \textit{e.g.} \cite[Chapter 3]{IwKo04}). 

A quasimodular form admits a Fourier expansion
\begin{equation}\label{eq:foudev}
f(z)=\sum_{n=0}^{+\infty}\widehat{f}(n)e(nz).
\end{equation}
Since we shall need to compute the odd part of a quasimodular form,
we introduce the notion of twist of a quasimodular form by a
Dirichlet character.

\begin{dfn}
Let $\chi$ be a Dirichlet character.  Let $f$ be a function having
Fourier expansion of the form \eqref{eq:foudev}.  The twist of $f$ by
$\chi$ is the function $f\otimes\chi$ defined by the Fourier expansion
\[
f\otimes\chi(z)
=
\sum_{n=0}^{+\infty}\chi(n)\widehat{f}(n)e(nz). 
\]
\end{dfn}

The interest of this definition is that it allows to build
quasimodular forms, as stated in the next proposition.
\begin{prop}\label{prop:twist}
Let $\chi$ be a Dirichlet character of conductor $m$ with non vanishing Gauss sum.  Let
$f$ be a quasimodular form of weight $k$ and depth $s$ on
$\Gamma_0(N)$.  Then $f\otimes\chi$ is a quasimodular form of weight
$k$, depth less than or equal $s$ and character $\chi^2$ on
$\Gamma_0\left(\lcm(N,m^2)\right)$.
\end{prop}
\begin{rem}
The Gauss sum of a character $\chi$ modulo $m$ is defined by
\[
\tau(\chi)=\sum_{u\ppmod{m}}\chi(u)e\left(\frac{u}{m}\right).
\]
\end{rem}
\begin{proof}
The proof is an adaptation of the corresponding result for modular
forms (see \textit{e.g.} \cite[Theorem 7.4]{Iwa97}).
We consider for each $k$ the following action of
$\mathrm{SL}_2(\R)$ on holomorphic functions on $\pk$:
\[
\trsf{f}{k}{\bigl(\begin{smallmatrix} a & b\\ c & d\end{smallmatrix}\bigr)}(z)
=
(cz+d)^{-k}f\left(\frac{az+b}{cz+d}\right).
\]
Since $\chi$ is primitive, this sum is not zero. One has
\begin{equation}
\label{eq:twistbar}
\tau(\overline{\chi})g\otimes\chi
=
\sum_{v\ppmod{m}}\overline{\chi}(v)\trsf{g}{k}{\bigl(\begin{smallmatrix} 1 & v/m\\ 0 & 1\end{smallmatrix}\bigr)}
\end{equation}
as soon as $g$ has a Fourier expansion of the form \eqref{eq:foudev}.
Define $M=\lcm(N,m^2).$  Let
$\bigl(\begin{smallmatrix} \alpha & \beta\\ \gamma &
\delta\end{smallmatrix}\bigr)\in\Gamma_0(M)$.  The matrix
\[
\begin{pmatrix}
1 & v/m\\ 0 & 1
\end{pmatrix}
\begin{pmatrix}
\alpha & \beta\\ \gamma & \delta
\end{pmatrix}
\begin{pmatrix}
1 & v\delta^2/m\\ 0 & 1
\end{pmatrix}^{-1}
\]
being in $\Gamma_0(N)$, one deduces from the level $N$ quasimodularity
of $f$ and \eqref{eq:twistbar} that
\begin{multline*}
\tau(\overline{\chi})\trsf{f\otimes\chi}{k}{\bigl(\begin{smallmatrix} \alpha & \beta\\ \gamma & \delta\end{smallmatrix}\bigr)}(z)
=\\
\sum_{i=0}^s\sum_{v\ppmod{m}}\overline{\chi}(v)
\trsf{f_i}{k-2i}{\bigl(\begin{smallmatrix} 1 & \delta^2v/m\\ 0 & 1\end{smallmatrix}\bigr)}(z)
\left[\frac{\gamma}{\gamma\left(z+\frac{\delta^2v}{m}\right)+\delta-\frac{\gamma\delta^2v}{m}}\right]^i.
\end{multline*}
Since the functions $f_i$ are themselves quasimodular forms (see \cite[Lemme 119]{MR05}), they admit a Fourier expansion.
Hence, from \eqref{eq:twistbar},
\[
\tau(\overline{\chi})\trsf{f\otimes\chi}{k}{\bigl(\begin{smallmatrix} \alpha & \beta\\ \gamma & \delta\end{smallmatrix}\bigr)}(z)
=
\tau(\overline{\chi})
\chi(\delta)^2
\sum_{i=0}^sf_i\otimes\chi(z)
\left(\frac{\gamma}{\gamma z+\delta}\right)^i.
\]
It follows that $f\otimes\chi$ satisfies the quasimodularity
condition.  There remains to prove the holomorphy at the cusps, which
is quite delicate since $f_s\otimes\chi$ may be $0$ even though $f_s$
is not.  Actually, lemma~\ref{lem:0040} and the fact that the twist of
a modular form on $\Gamma_0(N)$ by a primitive Dirichlet character of
conductor $m$ is a modular form on $\Gamma_0(M)$ show that the
proposition is proved as soon as it is proved for $f=D^{k/2-1}E_2$.
In that case, $s=k/2$ and $f_s\otimes\chi$ is not $0$ (see \cite[Lemma
118]{MR05}), hence $f_s$ being a modular form implies that
$f_s\otimes\chi$ is also one.

\end{proof}

\section{Proof of Hubert \& Leli\`evre conjecture}\label{sec:hlp}

The aim of this part is the proof of theorem~\ref{thm:pal}.  In all
this part, $n$ is assumed to be odd.  Define
\begin{align*}
\alpha_1(n,r)
&=
\sum_{(h_1,u_1,h_2,u_2)\in\Ar_1(n,r)}u_1u_2\\
\alpha_2(n,r)
&=
\frac{1}{2}\sum_{(h_1,u_1,h_2,u_2)\in\Ar_2(n,r)}u_1u_2\\
\alpha_3(n,r)
&=
\frac{n}{3}\sum_{(u_1,u_2,u_3)\in\Ar_3(n,r)}1
\end{align*}
with
\begin{multline*}
\Ar_1(n,r)
=\\
\left\{ (h_1,u_1,h_2,u_2)\in\Zp^4\colon
\begin{array}{|c}(h_1,h_2)=1,\\ \text{$h_1$ and $h_2$ odd,
}\end{array} u_1 < u_2,\, h_1u_1+h_2u_2=\frac{n}{r} \right\}
\end{multline*}
\begin{multline*}
\Ar_2(n,r)
=\\
\left\{
(h_1,u_1,h_2,u_2)\in\Zp^4\colon 
\begin{array}{|c}(h_1,h_2)=1,\\ \text{$h_1$ or $h_2$ even, }\end{array} 
\begin{array}{|c}u_1 < u_2,\\  \text{$u_1$ or $u_2$ even, }\end{array}
h_1u_1+h_2u_2=\frac{n}{r}
\right\}
\end{multline*}
and
\[
\Ar_3(n,r)
=
\left\{
(u_1,u_2,u_3)\in(2\Zpn+1)^3 \colon u_1+u_2+u_3=\frac{n}{r}
\right\}.
\]
By lemma \ref{lem:prph} and propositions \ref{prop:countUn} and \ref{prop:countDeux}, our goal is the
computation of
\begin{equation}\label{eq:but}
a_n^{\mathrm{p}}
=
\sum_{r\mid n}\mu(r)
\left[
r\alpha_1(n,r)+r\alpha_2(n,r)+\alpha_3(n,r)
\right].
\end{equation}
This, and hence theorem~\ref{thm:pal} is the consequence of the following lemmas \ref{lem:oddh}, \ref{lem:1948} and  \ref{lem:1949}.
\subsection{A preliminary arithmetical result}\label{sec:prel}
The following lemma will be useful in the sequel.
\begin{lem}\label{lem:1111}
Let $n\geq 1$, then
\[
\sum_{r\mid n}r\mu(r)
\sum_{d\mid n/r}\mu(d)\sigma_k\left(\frac{n}{rd}\right)
=
n^k
\sum_{r\mid n}
\frac{\mu(r)}{r^{k-1}}.
\]
and
\[
\sum_{r\mid n}\mu(r)
\sum_{d\mid n/r}\frac{\mu(d)}{d}\sigma_1\left(\frac{n}{rd}\right)
=
n\sum_{d\mid n}\frac{\mu(d)}{d^2}.
\]
\end{lem}
\begin{proof}
Consider the function
\[
f=(\id^\ell\mu)\ast\mu\ast\sigma_k.
\]
Then,
\[
L(f,s)=\frac{\zeta(s-k)}{\zeta(s-\ell)}=L(\id^\ell\mu,s)L(\id^k,s)
\]
hence
\[
f=(\id^\ell\mu)\ast\id^k.
\]
The lemma follows by taking $\ell=1$ for the first equality and $\ell=k=1$ for the second.
\end{proof}
\begin{rem}
Note that
\[
\sum_{r\mid n}\frac{\mu(r)}{r^2}
=
\prod_{p\mid n}\left(1-\frac{1}{p^2}\right).
\]
\end{rem}
\subsection{Two cylinders and odd heights}
Here, we compute the sum
\[
\sum_{r\mid n}r\mu(r)\alpha_1(n,r).
\]
More precisely, we prove the following lemma.
\begin{lem}\label{lem:oddh}
The number of type A primitive surfaces with $n$ squares and two cylinders of odd height is
\[
\frac{n^2(n-1)}{8}\sum_{r\mid n}\frac{\mu(r)}{r^2}.
\] 
\end{lem}
Write
\begin{equation}\label{eq:1333}
\alpha_1(n,r)
=
\gamma_1(n,r)-\widetilde{\alpha}_1(n,r)
\end{equation}
with
\[
\gamma_1(n,r)
=
\sum_{(h_1,u_1,h_2,u_2)\in\Cr(n,r)}u_1u_2
\]
and
\[
\widetilde{\alpha}_1(n,r)
=
\sum_{(h_1,u_1,h_2,u_2)\in\widetilde{\Ar}_1(n,r)}u_1u_2
\]
where
\[
\Cr(n,r)
=
\left\{
(h_1,u_1,h_2,u_2)\in\Zp^4
\colon
(h_1,h_2)=1,\,
u_1<u_2,\,
h_1u_1+h_2u_2=\frac{n}{r}
\right\}
\]
and (recalling that $n$ is odd)
\[
\widetilde{\Ar}_1(n,r)
=
\left\{
(h_1,u_1,h_2,u_2)\in\Zp^4\colon 
\begin{array}{|c}(h_1,h_2)=1,\\ \text{$h_1$ or $h_2$ even, }\end{array} 
u_1 < u_2,\, h_1u_1+h_2u_2=\frac{n}{r}
\right\}.
\]
Note that the sum
\[
\sum_{r\mid n}r\mu(r)\gamma_1(n,r)
\]
is the total number of primitive surfaces with two cylinders. Lemma \ref{lem:oddh} is a consequence of the two following lemmas
\ref{lem:two}, \ref{lem:evenh} and of equation \eqref{eq:1333}.
\subsubsection{Surfaces with two cylinders}
We prove the following result.
\begin{lem}\label{lem:two}
For $n$ odd, 
the number of primitive surfaces with $n$ squares and two cylinders is
\[
\frac{n^2(5n-18)}{24}\sum_{r\mid n}\frac{\mu(r)}{r^2}-\frac{n}{2}\varphi(n)
\]
where $\varphi$ is the Euler function.
\end{lem}

Using M{\"o}bius inversion formula, one obtains
\begin{align}
\notag
\gamma_1(n,r)
&=
\sum_{d\mid n/r}\mu(d)
\sum_{\substack{(i_1,u_1,i_2,u_2)\in\Zp^4\\ u_1<u_2\\ i_1u_1+i_2u_2=n/(rd)}}
u_1u_2\\
&=
\gamma_{1,1}(n,r)-\gamma_{1,2}(n,r)
\label{eq:1130}
\end{align}
with
\[
\gamma_{1,1}(n,r)
=
\frac{1}{2}\sum_{d\mid n/r}\mu(d)
\sum_{\substack{(i_1,u_2,i_2,u_2)\in\Zp^4\\ i_1u_1+i_2u_2=n/(rd)}}
u_1u_2
\]
and 
\[
\gamma_{1,2}(n,r)
=
\frac{1}{2}\sum_{d\mid n/r}\mu(d)
\sum_{\substack{(i_1,i_2,u)\in\Zp^3\\ (i_1+i_2)u=n/(rd)}}
u^2.
\]
One has
\[
\gamma_{1,1}(n,r)
=
\frac{1}{2}
\sum_{d\mid n/r}\mu(d)
\sum_{\substack{(v_1,v_2)\in\Zp^2\\v_1+v_2=n/(rd)}}
\sum_{w_1\mid v_1}w_1
\sum_{w_2\mid v_2}w_2
=
\frac{1}{2}
\sum_{d\mid n/r}\mu(d)
S_1\left(\frac{n}{rd}\right).
\]
By proposition~\ref{prop:essk}, this can be linearised to
\begin{multline*}
\gamma_{1,1}(n,r)
=
\frac{5}{24}\sum_{d\mid n/r}\mu(d)\sigma_3\left(\frac{n}{rd}\right)
-\frac{n}{4r}\sum_{d\mid n/r}\frac{\mu(d)}{d}\sigma_1\left(\frac{n}{rd}\right)
\\
+\frac{1}{24}\sum_{d\mid n/r}\mu(d)\sigma_1\left(\frac{n}{rd}\right)
\end{multline*}
so as to obtain
\begin{equation}\label{eq:1131}
\sum_{r\mid n}
r\mu(r)\gamma_{1,1}(n,r)
=
\left(\frac{5}{24}n^3-\frac{1}{4}n^2\right)
\sum_{r\mid n}\frac{\mu(r)}{r^2}
\end{equation}
thanks to lemma~\ref{lem:1111}.

Next, one has
\[
\gamma_{1,2}(n,r)
=
\frac{1}{2}
\sum_{d\mid n/r}
\mu(d)
\sum_{v\mid n/(rd)}
v^2
\left(\frac{n}{rdv}-1\right)
\]
so that
\begin{align}
\notag
\sum_{r\mid n}r\mu(r)\gamma_{1,2}(n,r)
&=
\frac{n}{2}\sum_{r\mid n}\mu(r)
\sum_{d\mid n/r}\frac{\mu(d)}{d}\sigma_1\left(\frac{n}{rd}\right)
\\
\notag
&\phantom{======}
-\frac{1}{2}\sum_{r\mid n}r\mu(r)\sum_{d\mid n/r}\mu(d)\sigma_2\left(\frac{n}{rd}\right)
\\
\label{eq:1132}
&=
\frac{n^2}{2}\sum_{r\mid n}\frac{\mu(r)}{r^2}
-
\frac{n}{2}\varphi(n)
\end{align}
by lemma~\ref{lem:1111}.

Finally, reporting \eqref{eq:1132} and \eqref{eq:1131} in \eqref{eq:1130} leads to
lemma \ref{lem:two}.
\subsubsection{Even product of heights}
Let us now compute the contribution of $\widetilde{\alpha}_1(n,r)$.
\begin{lem}\label{lem:evenh}
The number of type A primitive surfaces with $n$ squares and two cylinders, one having even height is
\[
\frac{n^2(2n-15)}{24}\sum_{r\mid n}\frac{\mu(r)}{r^2}+\frac{n}{2}\varphi(n).
\]
\end{lem}

Write
\begin{equation}\label{eq:1328}
\widetilde{\alpha}_1(n,r)
=
\widetilde{\alpha}_{1,1}(n,r)
-
\widetilde{\alpha}_{1,2}(n,r)
\end{equation}
with, recalling again that $n$ is odd,
\[
\widetilde{\alpha}_{1,1}(n,r)
=
\frac{1}{2}\sum_{d\mid n/r}\mu(d)\sum_{(i_1,u_1,i_2,u_2)\in\widetilde{\Ar}_{1,1}(n,r)}u_1u_2
\]
and
\[
\widetilde{\alpha}_{1,2}(n,r)
=
\frac{1}{2}\sum_{d\mid n/r}\mu(d)\sum_{(i_1,i_2,u)\in\widetilde{\Ar}_{1,2}(n,r)}u^2
\]
where
\[
\widetilde{\Ar}_{1,1}
=
\left\{
(i_1,u_1,i_2,u_2)\in\Zp^4 \colon \text{$i_1$ or $i_2$ even},\, i_1u_1+i_2u_2=\frac{n}{dr}
\right\}
\]
and 
\[
\widetilde{\Ar}_{1,2}
=
\left\{
(i_1,i_2,u)\in\Zp^3 \colon \text{$i_1$ or $i_2$ even},\, (i_1+i_2)u=\frac{n}{dr}
\right\}.
\]
Since $i_1$ and $i_2$ are not simultaneously even, one has
\[
\widetilde{\alpha}_{1,1}(n,r)
=
\sum_{d\mid n/r}\mu(d)\sum_{\substack{(v_1,v_2)\in\Zp^2\\ v_1+v_2=n/(dr)}}
\sum_{\substack{i_1\mid v_1\\\text{$i_1$ even}}}
\sum_{i_2\mid v_2}i_2
=
\sum_{d\mid n/r}\mu(d)S_2\left(\frac{n}{dr}\right).
\]
Using proposition~\ref{prop:essk} and lemma~\ref{lem:1111}, one obtains
\begin{equation}\label{eq:1329}
\sum_{r\mid n}r\mu(r)\widetilde{\alpha}_{1,1}(n,r)
=
\left(\frac{1}{12}n^3-\frac{1}{8}n^2\right)\sum_{r\mid n}\frac{\mu(r)}{r^2}.
\end{equation}
Next,
\[
\widetilde{\alpha}_{1,2}(n,r)
=
\frac{1}{2}\sum_{d\mid n/r}\mu(d)\sum_{u\mid n/(dr)}u^2\left(\frac{n}{rdu}-1\right).
\]
Lemma~\ref{lem:1111} gives
\begin{equation}\label{eq:1330}
\sum_{r\mid n}r\mu(r)\widetilde{\alpha}_{1,2}(n,r)
=
\frac{n^2}{2}\sum_{r\mid n}\frac{\mu(r)}{r^2}-\frac{n}{2}\varphi(n).
\end{equation}

Reporting \eqref{eq:1330} and \eqref{eq:1329} in \eqref{eq:1328} leads to
lemma \ref{lem:evenh}.

\subsection{Two cylinders with even product of heights}
Compute at last the sum
\[
\sum_{r\mid n}r\mu(r)\alpha_2(n,r).
\]
\begin{lem}\label{lem:1948}
The number of type A primitive surfaces with $n$ squares and two cylinders, one having an even height, the 
other having an even length is
\[
\frac{n^2(n-3)}{48}\sum_{r\mid n}\frac{\mu(r)}{r^2}.
\]
\end{lem}
Since $n$ is odd, one has
\[
\alpha_2(n,r)
=
\frac{1}{4}\sum_{d\mid n/r}\mu(d)
\sum_{(i_1,u_1,i_2,u_2)\in\widehat{\Ar}_{1,2}(n,r,d)}u_1u_2
\]
with
\[
\widehat{\Ar}_{1,2}(n,r,d)
=
\left\{
(i_1,u_1,i_2,u_2)\in\Zp^4\colon 
\begin{array}{|c}\text{$i_1$ and $u_1$ even},\\\text{or}\\ \text{$i_2$ and $u_2$ even, }\end{array} 
i_1u_1+i_2u_2=\frac{n}{rd}
\right\}.
\]
Hence,
\begin{align*}
\alpha_2(n,r)
&=
\frac{1}{2}
\sum_{d\mid n/r}\mu(d)\sum_{\substack{(i_1,u_1,i_2,u_2)\in\Zp^4\\\text{$i_1$ and $u_1$ even}\\i_1u_1+i_2u_2=n/(rd)}}u_1u_2\\
&=
\sum_{d\mid n/r}\mu(d)S_4\left(\frac{n}{rd}\right).
\end{align*}
The result follows from
proposition~\ref{prop:essk} and lemma~\ref{lem:1111}.
\subsection{One cylinder}
The counting in that case is more direct. 
\begin{lem}\label{lem:1949}
The number of type A primitive surfaces with $n$ squares and one cylinder is
\[
\frac{n^3}{24}\sum_{r\mid n}\frac{\mu(r)}{r^2}.
\]
\end{lem}
One actually has
\begin{align*}
\sum_{r\mid n}\mu(r)\alpha_3(n,r)
&=
\frac{n}{3}\sum_{r\mid n}\mu\left(\frac{n}{r}\right)
\#\left\{
(v_1,v_2,v_3)\in\Zpn^3 \colon v_1+v_2+v_3=\frac{r-3}{2}
\right\}
\\
&
=
\frac{n}{3}\sum_{r\mid n}\mu\left(\frac{n}{r}\right)\frac{r^2-1}{8}\\
&
=
\frac{1}{24}n^3\sum_{r\mid n}\frac{\mu(r)}{r^2}.
\end{align*}
\subsection{Computation of a generating series}

The number of non necessarily primitive surfaces with an odd number
$n$ of squares of type A is given by
\[
a_n
=
\sum_{d\mid n}\sigma_1\left(\frac{n}{d}\right)a_d^{\mathrm{p}}.
\]
Even though this does not have any geometric sense, one can define
numbers $a_n^\mathrm{p}$ and $a_n$ by these formulae for even
$n\geq2$.  We will compute the Fourier series attached to the
resulting sequence $\left(a_n\right)_{n\in\Zp}$.
Corollary~\ref{cor:qmodat} follows directly from the following
proposition.
\begin{prop}\label{prop:1421}
Let $n\geq 1$. Then
\[
a_n
=
\frac{3}{16}[\sigma_3(n)-n\sigma_1(n)].
\]
\end{prop}
\begin{proof}
We use the basic facts of \S~\ref{sec:prel}.
We have
\[
a_n=\frac{3}{16}[k_3(n)-k_2(n)]
\]
where, for $\ell\in\mathbb{Z}$, the arithmetical function $k_\ell$ is defined by
\[
k_\ell=\sigma_1\ast(\id^\ell\Psi)
\]
with
\[
\Psi(n)=\sum_{r\mid n}\frac{\mu(r)}{r^2}=\un\ast(\id\Psi)(n).
\]
We deduce that
\[
L(k_\ell,s)=\frac{\zeta(s)\zeta(s-1)\zeta(s-\ell)}{\zeta(s-\ell+2)}
\]
hence
\[
k_3=\sigma_3 \text{ and } k_2=\id\sigma_1.
\]
\end{proof}
\section{The associated Fourier series}\label{sec:labonneserie}
Recall that the two weight $2$ modular forms $\Phi_2$ and $\Phi_4$ on
$\Gamma_0(2)$ and $\Gamma_0(4)$ respectively have been defined in
\eqref{eq:lesphi} and \eqref{eq:losphi}.  In this section, we prove
theorem~\ref{thm:labonneserie}.
Since we want to eliminate the coefficients of even order, it is
natural to consider the Fourier series obtained by twisting all
coefficients by a modulus $2$ character.  By
proposition~\ref{prop:twist}, one obtains a quasimodular form of
weight $4$, depth less than or equal to $2$ on $\Gamma_0(4)$, hence a
linear combination of $E_4$, $E_{2,4}$, $E_{4,4}$, $D\Phi_2$,
$D\Phi_4$ and $DE_2$ (see the proof of proposition~\ref{prop:essk}).
The coefficients of this combination are found by computation of the
first seven Fourier coefficients.
\nocite{*}
\bibliographystyle{amsalpha}
\bibliography{BiblioLR05}
\end{document}